\documentclass[10pt,a4paper,twoside,fleqn,reqno]{rgn-art}
\usepackage{graphicx,enumerate}
\usepackage[dvips]{epsfig}
\usepackage[a4paper,textwidth=4.75in,textheight=7.75in,centering]{geometry}
\usepackage{hyperref} 

\theoremstyle{plain}
\newtheorem{thm}{Theorem}[section]

\newtheorem{lem}[thm]{Lemma}

\theoremstyle{definition}

\newtheorem{rem}[thm]{Remark}

\numberwithin{equation}{section}

\def\rrao{}
\begin{document}
\logo{
}{x (201x)}{x}

\title[Factorization of polynomials with analytic coefficients]{Factorization of polynomials with analytic coefficients}
\author{Wayne Lawton}

\address{\emph{Department of Mathematics, Mahidol University, Bangkok 10400, Thailand;
School of Mathematics and Statistics, University of Western Australia, Perth, Australia}.}
\email{scwlw@mahidol.ac.th}


\date{}

\begin{abstract}
We study monic univariate polynomials whose coefficients
are analytic functions of a real variable and whose roots lie in a specified
analytic curve. These include characteristic polynomials of unitary and
hermitian matrices whose entries are analytic functions. We use a result
of Newton to prove that every polynomial in such a class is a product of
degree one polynomials in the class.
\end{abstract}

\maketitle
\thispagestyle{empty}
\baselineskip13.5pt
\section{Introduction}
\label{sec1}
$\mathbb{R}$ and $\mathbb{C}$
are the real and complex numbers and
$\mathbb{T} = \{ \, z \in \mathbb{C} \, : \, |\, z\, | = 1 \, \}$ is the unit circle.
Functions defined by their Taylor series are called analytic.
For $r > 0,$ $A\, (\mathbb{D}_r)$ is the ring of analytic functions on the open disc
$\mathbb{D}_r =  \{ \, z \in \mathbb{C} \, : \, |z| < r \, \}$
and $A\, ((-r,r))$ is the ring of analytic functions on the open interval $(-r,r).$
We let $\mathbb{C}[z] \subset \mathcal{C}_{0}^{\omega} \subset \mathbb{C}[[z]]$ denote the rings of
polynomials, power series with complex coefficients that are absolutely convergent in $\mathbb{D}_r$ for some $r > 0,$ and formal power series. We identify $\mathcal{C}_{0}^{\omega}$ with the rings of germs of functions in $\cup_{\, r > 0} \ A\, ((-r,r))$ and of functions in $\cup_{\, r > 0} \ A\, (\mathbb{D}_r).$
\\
\indent $\mathcal{C}_{0}^{\omega}[z]$ is the ring of
polynomials with coefficients in $\mathcal{C}_{0}^{\omega}.$ Let $P(z) \in \mathcal{C}_{0}^{\omega}[z]$
be a monic polynomial of degree $d \geq 1.$ Then there exist $r > 0$ and
$a_0,...,a_{d-1} \in A\, (\mathbb{D}_r)$ such that $P(z) = z^d + a_{d-1}\, z^{d-1} + \cdots + a_1 \, z + a_0 \in \mathcal{C}_{0}^{\omega}[z].$ For $w \in \mathbb{D}_r$ we define $P_w(z) \in \mathbb{C}[z]$ by
$
P_w(z) = z^d + a_{d-1}(w)\, z^{d-1} + \cdots + a_1(w)\, z + a_0(w).
$
If $\gamma \subset \mathbb{C}$ we say that $P(z)$ {\it has roots in} $\gamma$ if there
exists $s \in (0,r]$ such that for every $t \in (-s,s)$ all roots of $P_t(z)$ are in $\gamma.$ We say that
$P(z)$ is {\it completely reducible} if factors
into monic polynomials in $\mathcal{C}_{0}^{\omega}[z]$ having degree one, or equivalently, if
there exist $u \in (0,s]$ and $\lambda_1,...,\lambda_d \in A\, (\mathbb{D}_u)$ such that
for every $w \in \mathbb{D}_u,$ $\lambda_1(w),....,\lambda_d(w)$ are the roots (with multiplicity) of $P_w(z).$ The polynomial $z^2 - t^2$ is completely reducible but the polynomial $z^2 - t$ is not.
In Section \ref{sec3} we prove:
\begin{thm}
\label{thm:main}
Every monic $P(z) \in \mathcal{C}_{0}^{\omega}[z]$
that has roots in an analytic curve $\gamma \subset \mathbb{C}$ is completely reducible.
\end{thm}
\section{Preliminary Results}
\label{sec2}
$\gamma \subset \mathbb{C}$ is an analytic curve if it is a real analytic submanifold of dimension $1.$ This means that for every point $p \in \gamma$ there exist $\epsilon > 0,$ an open neighborhood $U$ of $p$ in $\gamma,$ and an analytic diffeomorphism $f : (-\epsilon,\epsilon) \rightarrow U$ with $f(0) = p.$ Then $f^{\prime}(0) \neq 0.$
For $z = x + iy \in \mathbb{C}$ with $x, y \in \mathbb{R}$ we define $\Re \, z = x$ and $\Im \, z = y.$
\begin{lem}
\label{lem:imagreal}
If $\gamma \subset \mathbb{C}$ is an analytic curve and $p \in \gamma,$ then there exist
$c \in \mathbb{C} \, \backslash \, \{0\},$ $\delta > 0,$
an open neighborhood $V$ of $p$ in $\gamma,$
and an analytic function
$h \, : \, (-\delta,\delta) \rightarrow \mathbb{R}$
such that $h(0) = 0,$ $h^{\prime}(0) = 0,$ and
\begin{equation}
\label{eqn:imagreal}
 \Im \, \left[ \frac{z - p}{c} \right] = h \, \left( \, \Re \left[ \frac{z - p}{c}  \, \right] \, \right), \ \  \ \ z \in V.
\end{equation}
\end{lem}
\proof
Since $\gamma$ is an analytic curve and $p \in \gamma,$ there exist $\epsilon > 0,$ an open neighborhood $U$ of $p$ in $\gamma,$ and an analytic diffeomorphism $f : (-\epsilon,\epsilon) \rightarrow U$ such that
$f(0) = p.$ Let $c = f^{\prime}(0).$ Then $c \neq 0.$ Construct
$g : (-\epsilon, \epsilon) \rightarrow \mathbb{C}$
by $g(t) = (f(t)-p)/c$ and
$\psi = \Re \, g.$
Since $\psi(0) = 0$ and $\psi^{\prime}(0) = 1,$
the implicit function theorem for real analytic functions (\cite{2}, Theorem 1.4.3)
implies that there exist $\delta > 0$ and an analytic function
$\phi : (-\delta,\delta) \rightarrow (-\epsilon,\epsilon)$
such that $\phi(0) = 0,$ $\phi^{\prime}(0) = 1,$ and
$\psi(\phi(t)) = t, \, t \in (-\delta,\delta).$
Construct $h : (-\delta,\delta) \rightarrow \mathbb{R}$ by
$
h(t) = \Im \, g(\phi(t)).
$
Therefore
$h(0) = \Im \, g(\phi(0)) = \Im \, g(0) = \Im \, 0 = 0$
and
$h^{\prime}(0) = \Im \, (g^{\prime}(0) \, \phi^{\prime}(0)) = \Im \, 1 = 0.$
Let $V = f(\phi((-\delta,\delta))).$ Then $V$ is an open neighborhood of $p$ in $\gamma,$ and for every
$z \in V$ there exists $t \in (-\delta,\delta)$ with
$
z = f(\phi(t)).
$
Therefore
$$
\frac{z - p}{c} = \frac{f(\phi(t)) - p}{c} = g(\phi(t)).
$$
Equation \ref{eqn:imagreal} follows since
$
\Im \, g(\phi(t)) =
h(t) = h \left( \, \psi(\phi(t) \, \right) =
h \left( \, \Re \, g(\phi(t)) \right).
$
\endproof
\begin{lem}
\label{lem:newton}
If $P(z)$ is a monic polynomial that is irreducible in $\mathcal{C}_{0}^{\omega}[z]$ and
has degree $d \geq 2$ then there exist $r > 0$ and $\eta \in A\, (\mathbb{D}_r)$ such that
\begin{equation}
\label{eqn:fac1}
    P_{w^d}(z) = \prod_{k=0}^{d-1} \left[\, z - \eta(e^{2\, \pi\, i\, k/d}\, w)\, \right], \ \ w \in \mathbb{D}_r.
\end{equation}
\end{lem}
\proof
Abhyankar (\cite{2}, Newton's Theorem and Supplements 1 and 2 on page 89) proves a version of this result for polynomials with coefficients in the ring of formal power series $\mathbb{C}[[w]]$ and says that it was proved by Newton in 1660 \cite{5}. The version in Lemma \ref{lem:newton} for coefficients in
$\mathcal{C}_{0}^{\omega}$ follows from Weierstrass' $M$-test.
\endproof
\begin{lem}
\label{lem:eta}
If $\eta$ in Equation \ref{eqn:fac1} has the Taylor expansion
$\eta(w)= \sum_{\, n = 0}^{\, \infty}  \eta_n \, w^n,$
then there exists $L \geq 1$ such that $\eta_L \neq 0$ and $d$ does not divide $L.$
\end{lem}
\proof
Otherwise there exists $\mu \in A\, (\mathbb{D}_{r^d})$ such that
$\eta(w) = \mu(w^d), \ \ w \in \mathbb{D}_r.$
Then Equation \ref{eqn:fac1} implies that
$
P_{w^d}(z) = \left(\, z - \mu(w^d)\, \right)^d, \ \ w \in \mathbb{D}_r.
$
Since the function $w \rightarrow w^d$ maps $\mathbb{D}_r$ onto $\mathbb{D}_{r^d},$
$
P_{w}(z) = \left(\, z - \mu(w)\, \right)^d, \ \ w \in \mathbb{D}_{r^d},
$
so $P(z)$ is not irreducible in $\mathcal{C}_{0}^{\omega}[z].$ This contradiction completes the proof.
\endproof
\section{Proof of Theorem \ref{thm:main}}
\label{sec3}
Assume to the contrary that there exist an analytic curve $\gamma \subset \mathbb{C}$ and a monic polynomial $P(z) \in \mathcal{C}_{0}^{\omega}[z]$ of degree $d \geq 2$ that has roots in $\gamma$ and is not completely reducible. We may assume that $P(z)$ is irreducible in $\mathcal{C}_{0}^{\omega}[z]$ so Lemma \ref{lem:newton} implies there exist $r > 0$ and $\eta \in A\, (\mathbb{D}_r)$ that satisfy Equation \ref{eqn:fac1}. Since the roots of $P(z)$ are in $\gamma,$ there exists $s \in (0,r]$ such that
$\eta(w) \in \gamma$ whenever $w^d \in \mathbb{R}$ and $w \in \mathbb{D}_s.$
Let $p = \eta(0).$ Lemma \ref{lem:imagreal} implies that there exist
$c \in \mathbb{C} \, \backslash \, \{0\},$
$\delta > 0,$ an open neighborhood $V$ of $p$ in $\gamma,$ and an analytic function
$h \, : \, (-\delta,\delta) \rightarrow \mathbb{R}$
such that $h(0) = 0,$ $h^{\prime}(0) = 0,$ and Equation \ref{eqn:imagreal} holds.
Since $\eta$ is continuous there exists $u \in (0,s]$ such that
$\eta(w) \in V$ whenever $w^d \in \mathbb{R}$ and $w \in \mathbb{D}_u.$
Construct $\lambda = (\eta - p)/c$ with Taylor series
$\sum_{\, n = 0}^{\, \infty} \lambda_n \, w^n.$
Then $\lambda_0 = 0$ and Lemma \ref{lem:eta} implies that there exists a smallest integer
$L \geq 1$ such that $\lambda_L \neq 0$ and $d$ does not divide $L.$
Choose $k \in \{ \, 0,1,2,...,d-1\, \}$ such that
$\Im \, (e^{\pi i k L/d} \, \lambda_L) \neq 0$
and construct
$\zeta(t) = \Im \, \lambda(e^{\pi i k/d} t), \, t \in (-u,u)$
with Taylor series
$\sum_{\, n = 0}^{\, \infty} \zeta_n \, t^n.$
Then
$\zeta_L = \Im \, (e^{\pi i k L/d} \, \lambda_L) \neq 0.$
If $t \in (-u,u)$ then $\eta(e^{\pi i k L/d} t) \in V$ so Equation \ref{eqn:imagreal} gives
$
    \zeta(t) = h \, \left( \, \Re \, \lambda(e^{\pi i k L/d} t) \, \right).
$
The facts that $1 \leq m < L$ implies that $d$ divides $m$ or $\lambda_m = 0,$ $\lambda_0 = 0,$ $h^{\prime}(0) = 0,$ and $d$ does not divide $L,$ imply that $\zeta_L = 0.$
This contradiction completes the proof.
\begin{rem}
In (\cite{3}, Corollary 1) we proved that a monic $P(z) \in \mathcal{C}_{0}^{\omega}[z]$
of degree $2$ that has roots in $\mathbb{T}$ is completely reducible and used results in \cite{4} to prove that the eigenvalues of certain unitary matrices (arising in quantum physics) with analytic entries are global analytic functions on $\mathbb{T}$ if the characteristic polynomials of the matrices are completely reducible. Theorem \ref{thm:main} ensures this condition holds.
\end{rem}


\end{document}